# Modular analogues of Jordan's theorem for finite linear groups

Michael J. Collins

**Abstract.** In 1878, Jordan [**9**] showed that a finite subgroup of $GL(n, \mathbb{C})$ contains an abelian normal subgroup whose index is bounded by a function of $n$ alone. Previously, the author has given precise bounds [**4**]. Here, we consider analogues for finite linear groups over algebraically closed fields of positive characteristic $\ell$. A larger normal subgroup must be taken, to eliminate unipotent subgroups and groups of Lie type and characteristic $\ell$, and we show that generically the bound is similar to that in characteristic 0 - being $(n+1)!$, or $(n+2)!$ when $\ell$ divides $n+2$ - given by the faithful representations of minimal degree of the symmetric groups. A complete answer for the optimal bounds is given for all degrees $n$ and every characteristic $\ell$.

MSC (2000): 20C20, 20C33, 20D25

## 1. Introduction

A theorem due to Jordan [**9**] states that there is a function $f$ on the natural numbers such that, if $G$ is a finite subgroup of $GL(n, \mathbb{C})$, then $G$ has an abelian normal subgroup $N$ of index bounded by $f(n)$. In two previous papers ([**3**],[**4**]), we have given explicit bounds and have shown that they can be achieved; in particular, the generic bound is $(n+1)!$, achieved by the symmetric group, and this holds whenever $n \geq 71$.

In this paper, we will study the corresponding situation for finite subgroups of $GL(n, k)$ where, as throughout this paper, $k$ will denote an algebraically closed field of positive characteristic $\ell$. We have to take account of a number of major differences; we must allow for an arbitrarily large normal unipotent subgroup, and we must allow for groups of Lie type in characteristic $\ell$ over arbitrarily large finite fields. But there is a further consideration: if $\ell$ divides $m$, then the standard permutation module of the symmetric group $S_m$ in characteristic $\ell$ is uniserial with trivial head and socle, with a middle composition factor of dimension $m-2$. In this case, for $m \geq 5$, this is the smallest faithful representation.

This has two consequences, even after taking into account unipotent and characteristic $\ell$ groups. Clearly the bound $(n+1)!$ on some suitable quotient no longer holds uniformly for large enough $n$, specifically when $\ell$ divides $n+2$. But there is also a more subtle point: when $\ell$ divides $n+1$, we may (and will for large enough $n$) get a bound of $(n+1)!$, but the group $S_{n+1}$ is no longer irreducible and hence certainly not primitive.

This influences our strategy. Rather than determine the best possible bounds immediately, first for primitive groups and then for all groups as we did in characteristic zero, we will view $(n+1)!$ as the "generic" bound for the order of a particular quotient of $G$ that we will describe, and regard the situation when $\ell$ divides $n+2$ and the bound is $(n+2)!$ as an exception, and for each characteristic $\ell$ we will determine the smallest degree above which one of these must hold. Only then, in Section 10, and given the techniques that we adopted in our previous papers and sketch here, do we give precise bounds for smaller degrees (Theorems C – F). These results clearly depend on $\ell$; more seriously, had we attempted to integrate this with the earlier work, then we would have obscured our key aim of seeking the generic bound. Furthermore, to obtain the best bounds for *primitive* groups when $\ell$ divides $n+1$ would not be straightforward, while $(n+1)!$ turns out to be a perfectly adequate working bound.

Our main theorem is to obtain the following generic result, as an analogue of our earlier result in characteristic 0.

**Theorem A.** *Let $k$ be an algebraically closed field of positive characteristic $\ell$ and let $G$ be a finite subgroup of $GL(n,k)$. Put $\bar{G} = G / O_\ell(G)$. If $n \geq 71$, then $\bar{G}$ has normal subgroups $\bar{N}$ and $\bar{L}$ such that*

*(i)   $\bar{N}$ is abelian,*

*(ii)  $\bar{L} = E_\ell(\bar{G})$, and*

*(iii) $[\bar{G} : \bar{L}\bar{N}] \leq \begin{cases} (n+2)! & \text{if } \ell \text{ divides } n+2 \\ (n+1)! & \text{otherwise.} \end{cases}$*

*If the bound is achieved, then $E(\bar{G})$ is simple and $\bar{G} / Z(\bar{G}) \cong S_{n+2}$ or $S_{n+1}$, respectively. Furthermore, $O_\ell(G) = 1$ unless $n \equiv -1 \pmod{\ell}$, in which case $O_\ell(G)$ may be elementary abelian.*

We will define the subgroup $E_\ell(X)$ of an arbitrary finite group $X$ shortly, but first we remark that Theorem A has been stated in a form that is optimal but independent of the characteristic of the field $k$. If $k$ has small characteristic, we can achieve the generic bound in smaller degree.

**Theorem A′.** *If $\ell \leq 5$, then the bounds and subsequent conclusion of Theorem A hold whenever $n \geq n(\ell)$, where*

$$n(2) = 34, \ \ n(3) = 69, \ and \ n(5) = 70.$$

*For $\ell \geq 7$, the condition that $n \geq 71$ is optimal.*

We shall now introduce some standard terminology and notation. Let $X$ be any finite group and $p$ a prime. Then $O_p(X)$ is the largest normal $p$-subgroup of $X$. By a *component* of $X$ we mean a subnormal subgroup that is a perfect central extension of a nonabelian simple group (i.e., *quasisimple*), and the *Bender subgroup* $E(X)$ is the (central) product of all the components of $X$. Now let $\mathcal{L}_p$ denote the family of all finite simple groups of Lie type in characteristic $p$; then we put

$$E_p(X) = \langle E \lhd E(X) \mid E = E', \, E \, / \, Z(E) \in \mathcal{L}_p \rangle$$

and

$$E_{p'}(X) = \langle E \lhd E(X) \mid E = E', \, E \, / \, Z(E) \notin \mathcal{L}_p \rangle.$$

In particular, $E_p(X)$ and $E_{p'}(X)$ are characteristic subgroups of $E(X)$ and hence also of $X$.

In making this final definition, those simple groups that occur in more than one guise recur – in particular, we include the alternating group $A_5$ as a member of both $\mathcal{L}_2$ and $\mathcal{L}_5$. Absorbing $A_5$ into $\mathcal{L}_2$ is what makes $n(2)$ as small as it is in Theorem A′; otherwise we would have had $n(2) = 69$, as for $\ell = 3$ and almost as in characteristic zero.

Before his disappearance in 1985, Boris Weisfeiler obtained the bounds

$$[\bar{G} : \bar{L}\bar{N}] \leq \begin{cases} (n+2)! & \text{if } n \geq 64 \\ n^4(n+2)! & \text{if } n < 64 \end{cases}$$

in a near complete unpublished manuscript, covering all possible characteristics $\ell$ including zero. Our methods, which are rather different, describe the bounds precisely. It is then straightforward but tedious to verify his bound for small $n$ and, in the final section of this paper, we will consider the particular case where the extra factor $n^4$ turns out to be closest to optimal as an illustration, though in general it gives a poor "correction term". Separately, and using methods from algebraic geometry and linear algebraic groups but not requiring the classification of finite simple groups, Larsen and Pink have established the existence of bounds. (See [**10**].) We mention (as did Weisfeiler) that the bounds obtained by Brauer and Feit [**2**] follow from Theorem A.

We refer the reader to Section 47 of [**1**] for our notation for the finite simple groups, and to Section 31 for properties of $E(G)$ and related subgroups.

*Acknowledgements.* Much of this paper was written while the author was a Visiting Associate at the California Institute of Technology during the Winter term 2005. The author would like to thank the referee for a number of detailed comments that have been incorporated into this revision.

## 2. Reduction to the semisimple case

We start with a simplification.

**Proposition 1.** *Let $G$ be a finite subgroup of $GL(n,k)$ and put $\bar{G} = G / O_\ell(G)$. Then $\bar{G}$ has a faithful, completely reducible, representation of degree $n$.*

*Proof.* Let $V$ be the underlying vector space and let

$$V = V_1 \supset \cdots \supset V_{m+1} = 0$$

be a composition series for $V$ as a $kG$-module. Then the stabiliser of the chain is an $\ell$-subgroup of $G$ while $O_\ell(G)$ lies in the kernel of every simple $kG$-module. Thus $\bar{G}$ acts faithfully on the direct sum

$$\oplus_{i=1}^m V_i\,/\,V_{i+1}\,,$$

and this has dimension $n$. $\square$

Since the bounds in Theorem A make claims only about the quotient $\bar{G} = G\,/\,O_\ell(G)$, we may suppose throughout, without loss, that $O_\ell(G) = 1$ and that $G$ acts completely reducibly (except for Proposition 26). This will enable us to apply some of the same techniques as in the characteristic zero case, once we have obtained bounds for primitive groups in Theorem B in the next section and can take care of components of Lie type in characteristic $\ell$.

## 3. Bounds for primitive groups

Recall that an irreducible (linear) group is *primitive* if the underlying vector space does not decompose into a direct sum of proper subspaces permuted under the action of the group. This is equivalent to saying that the associated representation is not induced from any proper subgroup. By Clifford's theorem, a primitive linear group is also *quasiprimitive* – i.e., irreducible, with every normal subgroup acting homogeneously. In Section 2 of [**3**], we considered primitive groups in arbitrary characteristic; Proposition 2 collects some of those results.

**Proposition 2.** *Let $G$ be a primitive subgroup of $GL(n,k)$. Then the following hold.*

*(i) Every abelian normal subgroup of $G$ is cyclic and central.*

*(ii) Let $p$ be a prime different from $\ell$ and put $P = O_p(G)$. Then either $P$ is cyclic, or else $P$ contains an extraspecial subgroup $E$ such that $P = Z(P).E$.*

*Furthermore, if $P$ is noncyclic and $|\,E\,| = p^{2m+1}$, then*

*(iii) the stabiliser $C_G(P/Z(P)) \cap C_G(Z(P))$ of the chain $P \supset Z(P) \supset 1$ is just $P.C_G(P)$,*

*(iv) $G/P.C_G(P)$ is isomorphic to a subgroup of $Sp_{2m}(p)$, and*

*(v) if $p$ is odd, then $E$ has exponent $p$ and $E = \Omega_1(P)$  ( $= \langle x \in P \mid x^p = 1 \rangle$ ).*

We recall that the *generalised Fitting subgroup* of $G$ is

$$F^*(G) = F(G).E(G)$$

where $F(G)$ is the Fitting subgroup. A noncyclic group $P$ as in (ii) will be called a *quasicomponent*, and we let[1] $E^*(G)$ be the product of the components and quasicomponents of $G$. In particular, here we have $F^*(G) = Z(G).E^*(G)$ and, by Proposition 2 and the general property of the generalised Fitting subgroup that $C_G(F^*(G)) \subseteq F^*(G)$, that

$$C_G(E^*(G)) = Z(G).$$

We will require a refinement of Theorem 5 of [**3**]. Let $G$ be a primitive group with components $E_1, \ldots, E_s$. Suppose that $s \geq 1$. Then $E(G)$ acts homogeneously on the underlying vector space $V$. Let $U$ be a simple submodule of $V$, regarded as a $kE(G)$-module. Then $E(G)$ acts faithfully on $U$ and $U$ decomposes as a tensor product

$$U = U_1 \otimes \cdots \otimes U_s$$

of nonzero spaces where $E_j$ acts irreducibly on $U_j$ (not necessarily faithfully, but with kernel contained in $Z(E_j)$), and trivially on $U_{j'}$ for $j' \neq j$.

Put

$$N = \bigcap_{j=1}^{s} N_G(E_j).$$

For each $j$, $E_j \lhd N$ and $Z(E_j) \subseteq Z(G) = Z(N)$; hence there is a homomorphism $\varphi_j : N \to Out_c(E_j)$, the subgroup of $Out(E_j)$ which is the image of the group $Aut_c(E_j)$ of automorphisms of $E_j$ that act trivially on $Z(E_j)$, and the kernel of this homomorphism contains $F^*(G)$. Furthermore, since $E(G)$ acts homogeneously, $N$ must stabilise the representation $\rho_j$ of $E_j$ afforded by $U_j$ so that

$$\varphi_j(N) \subseteq I_{Out_c(E_j)}(\rho_j),$$

the inertia group of $\rho_j$ in $Out_c(E_j)$. Thus, together with Theorem 5 of [**3**], we have established the following, the refinement being the replacement of subgroups $Out_c(E_j)$ by the inertia groups $I_j$ in (ii).

---

[1]  We use $E^*(G)$ in this paper instead of $E_1(G)$ as in [**3**] and [**4**] to avoid confusion with $E_c(G)$.

**Theorem 3.** *Let $G$ be a nonabelian primitive group with quasicomponents $P_1,\dots,P_r$ and components $E_1,\dots,E_s$. For each $i$, put $|P_i/Z(P_i)| = p_i^{2m_i}$ and let $N = \bigcap_{j=1}^{s} N_G(E_j)$. For each $j$, let $I_j = I_{Out_c(E_j)}(\rho_j)$ where $\rho_j$ is defined above. Then*

*(i)    there is a monomorphism from $N/F^*(G)$ into the direct product*

$$Sp_{2m_1}(p_1) \times \cdots \times Sp_{2m_r}(p_r) \times I_1 \times \cdots \times I_s,$$

*and*

*(ii)   $G/N$ is isomorphic to a subgroup of a direct product $S_{l_1} \times \cdots \times S_{l_t}$ of symmetric groups where $l_1,\dots,l_t$ are the sizes of the distinct isomorphism classes of components of $G$.*

Bounds for primitive groups will be given by the following theorem; when $n > 12$ they are best possible except when $\ell$ divides $(n+1)$, and the explicitly listed cases for $n \le 12$ are also all optimal.

**Theorem B.** *Let $G$ be a finite primitive subgroup of $GL(n,k)$. Then*

$$[G : Z(G).E_\ell(G)] \le (n+1)!$$

*with the following exceptions.*

(i)  $\ell$ divides $(n+2)$ and $[G : Z(G).E_\ell(G)] \leq (n+2)!$, for $n \geq 9$, $n \neq 12$, or $n = 8$ when $\ell = 2$.

(ii)  $n \leq 9$ or $n = 12$ for appropriate $\ell$, where the bounds are given below, together with particular groups having centre of minimal order that achieve them:

| $n$ | $[G : Z(G).E_\ell(G)]$ | restriction on $\ell$ | |
|---|---|---|---|
| 2 | 60 | $\ell \neq 2,5$ | $2 \cdot A_5 \ (\cong 2 \cdot SL_2(4),\ SL_2(5))$ |
| | 24 | $\ell = 5$ | $2 \cdot S_4$ |
| 3 | 360 | $\ell \neq 2,3,5$ | $3 \cdot A_6 \ (\cong 3 \cdot Sp_4(2)',\ 3 \cdot L_2(9))$ |
| | 216 | $\ell = 2$ | $3^{1+2}.SL_2(3)$ |
| | 168 | $\ell = 3$ | $L_2(7) \ (\cong L_3(2))$ |
| | 2520 | $\ell = 5$ | $3 \cdot A_7$ |
| 4 | 25920 | $\ell \neq 2,3$ | $Sp_4(3) \ (\cong 2 \cdot U_4(2))$ |
| | 2520 | $\ell = 2$ | $A_7$ |
| | 40320 | $\ell = 3$ | $4_2 \cdot (L_3(4).2_2)$ |
| 5 | 25920 | $\ell \neq 2,3$ | $PSp_4(3) \ (\cong U_4(2))$ |
| | 3000 | $\ell = 2$ | $5^{1+2}.SL_2(5)$ |
| | 7920 | $\ell = 3$ | $M_{11}$ |
| 6 | 6531840 | $\ell \neq 2,3$ | $6_1 \cdot (U_4(3).2_2)$ |
| | 6531840 | $\ell = 2$ | $3_1 \cdot (U_4(3).2_2)$ |
| | 604800 | $\ell = 3$ | $2 \cdot J_2$ |
| 7 | 1451520 | $\ell \neq 2$ | $Sp_6(2)$ |
| 8 | 348368800 | $\ell \neq 2$ | $2 \cdot (O_8^+(2).2)$ |
| 9 | 4199040 | $\ell \neq 2,3,5$ | $3^{1+4}.Sp_4(3)$ |
| | 50232960 | $\ell = 2$ | $3 \cdot J_3$ |
| | 12700800 | $\ell = 5$ | $(3 \cdot A_7 \circ 3 \cdot A_7)Z_2$ |
| 12 | 448345497600 | $\ell \neq 2,3$ | $6 \cdot Suz$ |
| | 448345497600 | $\ell = 2$ | $3 \cdot Suz$ |
| | 896690995200 | $\ell = 3$ | $2 \cdot (Suz.2)$ |

*where the group in the final column is described by its the normal structure using the notation[2] of the ATLAS [5].*

*Furthermore, when the bound in any of the exceptional cases is achieved, $O_\ell(G) = E_\ell(G) = 1$.*

**Remarks.** (i) If $\ell \geq 7$, the bounds in the exceptional cases are all the same as in characteristic zero. Thus our basic argument will be the same, except that we need to take account of symmetric groups $S_{n+2}$ in degree $n$. Because the group $2 \cdot A_5$ appears in degree 2 for $\ell = 3$, the arguments there are similar, but the cases $\ell = 2$ and $\ell = 5$ will differ substantially.

(ii) The "generic" bound of $(n+1)!$ is achieved for $n > 12$ provided that $\ell$ does not divide either $n+1$ or $n+2$. When $\ell$ divides $n+1$, the bound $(n+1)!$ arises by reduction to a configuration that does not occur in a primitive group; the optimal bounds are not known.

(iii) In degrees 8 and 9, the bound is given by $S_{10}$ if $\ell = 2$ or 3, respectively. Degree 7 for $\ell = 2$ exhibits the central problem posed above in (ii); however, here we even have $A_8 \in \mathcal{L}_2$. The optimal bound in fact is 16464, given by the group $7^{1+2}.SL_2(7)$. We note also that the Mathieu group $M_{24}$ has a representation of degree 11 in characteristic 2 only and achieves the optimal bound in that case.

(iv) In degree 2 for $\ell = 2$, the bound is actually 1; every primitive group is of the form $Z(G) \times SL(2, 2^m)$.

(v) Since the groups listed in the final column will have a vital role to play in the general problem, we have described those that belong to $\mathcal{L}_p$ for more than one prime $p$.

---

[2] We have made one small change. For an arbitrary group $X$, we denote a nonsplit central extension by a cyclic subgroup of order $m$ by $m \cdot X$. Whenever we write $XY$ or $X.Y$, we will mean either a product of subgroups or a split extension of abstract groups (which will be clear from the context). Also, $(3 \cdot A_7 \circ 3 \cdot A_7)Z_2$ denotes the nondirect central product of two copies of $3 \cdot A_7$ interchanged by an involution.

Let $G$ be a primitive subgroup of $GL(n,k)$. We will first show that $[G : Z(G).E_\ell(G)]$ is bounded by a function of $n$. Since we wish to determine when the generic bound $(n+1)!$ is exceeded, we will suppose that $G$ is such that the index is maximal and

$$[G : Z(G).E_\ell(G)] > (n+1)!;$$

what we will then show is that either $\ell$ divides $(n+2)$ and $G/Z(G) \cong S_{n+2}$ or else one of the exceptions in (ii) occurs with $G/Z(G) \cong H/Z(H)$ where $H$ is the group listed. First, we note that all the groups listed do have primitive representations of the degrees claimed since an irreducible group that is not primitive necessarily has a suitable permutation group as a homomorphic image; except for $2 \cdot S_4$, $5^{1+2}.SL_2(5)$, $3^{1+4}.Sp_4(3)$ and $(3 \cdot A_7 \circ 3 \cdot A_7)Z_2$, the representations can be found in either the modular atlas [8] or are reductions $\bmod \ell$ of ordinary representations that can be found in the ATLAS [5], and primitivity is a consequence of the near simplicity of the groups involved. The group $2 \cdot S_4$ has $2 \cdot A_4$ as a subgroup of index 2, and its embedding in $SL_2(5)$ yields primitivity, while $3^{1+2}.SL_2(3)$, $5^{1+2}.SL_2(5)$ and $3^{1+4}.Sp_4(3)$ were shown in Proposition 7 of [3] to have primitive complex representations of degrees 3, 5 and 9 respectively. The modular representations claimed are then obtained by reduction $\bmod \ell$ and remain primitive. The group $3 \cdot A_7$ has an irreducible representation of degree 3 when $\ell = 5$, and the direct product of two copies of the group has an irreducible representation of degree 9 given by their tensor product which is necessarily primitive and extends naturally to the wreath product $(3 \cdot A_7)wrS_2$ with a nontrivial kernel.

Our proof of Theorem B will be modelled on that for the corresponding result in characteristic zero; however, first we need to study irreducible representations for extensions of groups of Lie type in characteristic $\ell$.

## 4. Extensions of linear quasisimple groups of Lie type in characteristic $\ell$

Let $G$ be a group for which $E(G)/Z(E(G)) \in \mathcal{L}_\ell$ and $F^*(G) = Z(G).E(G)$. Suppose that $G$ has an irreducible representation $\varphi$ over $k$ whose restriction to $E(G)$ is nontrivial and homogeneous. Since the exceptional Schur multipliers for groups of Lie type occur only in the defining characteristic and therefore lie in the kernel of any irreducible representation in that characteristic, we may assume without loss that $E(G)$ is the universal covering group of a simple group (in the sense of Steinberg [**12**]).

Put $E = E(G)$. Let $\rho$ be an irreducible constituent of $\varphi \mid_E$. Suppose that $\rho$ has degree $n$. We will bound the index $[G : F*(G)]$ in terms of $n$ (and hence in terms of the degree of $\varphi$). As in Section 3, conjugation in $G$ induces a map

$$G \to Out(E)$$

with kernel $F*(G)$, and $Out(E)$ permutes the irreducible representations of $E$; hence

$$[G : F*(G)] \leq \mid I_{Out(E)}(\rho) \mid.$$

Suppose that $E = E(q)$ where $q = \ell^b$ and $E(q)$ is a group of Lie type over $\mathbb{F}_q$ (possibly twisted, in which case recall that our notation differs from that in [**13**]). Let $\sigma : t \mapsto t^\ell$ be the Frobenius automorphism of $k$. Then $\sigma$ induces an automorphism of $E$ that we also denote by $\sigma$. Let $F = \langle \sigma \rangle$ be the group of field automorphisms of $E$; then $F \cong Z_{ab}$, where $a = 3$ for $E$ of type $^3D_4(q)$, $a = 2$ for types $^2A_r(q)$, $^2D_r(q)$ or $^2E_6(q)$, and $a = 1$ otherwise. $Out(E)$ has a normal series

$$Out(E) \supseteq FD \supseteq D \supseteq 1$$

where $D$ is the image in $Out(E)$ of the diagonal automorphisms, $F$ is viewed as a subgroup of $Out(E)$, and $Out(E) / FD$ is the group of graph automorphisms.

By Steinberg's tensor product theorem ([**13**], Theorem 13.3), there is a set of fundamental representations $P = \{\rho_\alpha \mid \alpha \in A\}$ such that every irreducible representation of $E$ is uniquely of the form

$$\otimes_{i=0}^{b-1} \rho_i^{\sigma^i}$$

where, letting $k$ denote the trivial representation, for each $i$ we have $\rho_i \in P \cup k$ and $\rho_i^{\sigma^i}(x) = \sigma^i(\rho_i(x))$. Furthermore, since each such representation is obtained by restriction from a representation of the ambient algebraic group $E(k)$, the diagonal automorphisms of $E$ stabilise every irreducible representation ([**11**], p.607); in particular $D \subseteq I_{Out(E)}(\rho)$, so that $I_{FD}(\rho) = I_F(\rho).D$.

The uniqueness of this tensor product decomposition yields an action of $F$ on the set of irreducible representations of $E$ via

$$\sigma : \otimes_{i=0}^{b-1} \rho_i^{\sigma^i} \mapsto \otimes_{i=0}^{b-1} \rho_i^{\sigma^{i+1}}$$

with $\rho_{b-1}^{\sigma^b} \in \mathrm{P} \cup k$. Thus $I_F(\rho)$ permutes the set of tensor factors $\rho_i^{\sigma^i}$ into cycles for each of which the fundamental representations involved are the same. We first consider the case where $a = 1$, i.e., $E(q)$ is untwisted or of type $^2B_2$, $^2F_4$ or $^2G_2$. Then, taking any cycle as described above consisting of images of a fundamental or twisted representation of degree $m$ and by considering the contribution to the degree of $\rho$, we see that $m^{|I_F(\rho)|}$ divides $n$.

When $a \neq 1$, similar considerations apply. However, now $F$ cycles through the tensor factors $a$ times and we can conclude only that $m^{|I_F(\rho)|/a}$ divides $n$.

**Proposition 4.** *Let $E = E(q)$ be the universal covering group of a finite simple group of Lie type of characteristic $\ell$. Let $\rho$ be an irreducible representation of $E$ in characteristic $\ell$ of degree $n$ and suppose that, in the tensor decomposition of $\rho$, the maximal degree of any fundamental representation occurring in the tensor decomposition is $m$. Then*

$$| I_{Out(\bar{E})}(\rho) | \leq \lambda a n \log_m n$$

*where*

$$\lambda = \begin{cases} 1 \text{ if } E \text{ is twisted} \\ 2 \text{ if } E \text{ is untwisted and not of type } D_4 \\ 6 \text{ if } E \text{ is of type } D_4. \end{cases}$$

*Proof.* In the notation established earlier, we have

$$| I_{Out(\bar{E})}(\rho) | \leq | D | \cdot | I_F(\rho) | \cdot [Out(E) : FD].$$

If $E$ is a classical group, $| D |$ divides either 4 or the natural degree; in any case, $| D | \leq n$. For the remaining types, $| D | \leq 3$ and the same conclusion holds. We have shown above that

$$| I_F(\rho) | \leq a \log_m n$$

while

$$[Out(\bar{E}) : FD] \leq \begin{cases} 1 \text{ if } E \text{ is twisted} \\ 2 \text{ if } E \text{ is untwisted and not of type } D_4 \\ 6 \text{ if } E \text{ is of type } D_4. \quad \square \end{cases}$$

**Corollary 5.** *Let $G$ be a group and assume that $F^*(G) = Z(G).E(G)$ where $E(G) / Z(E(G)) \in \mathcal{L}_\ell$. Suppose that $G$ has a faithful irreducible representation $\varphi$ over $k$ whose restriction to $E(G)$ is homogeneous. If an irreducible constituent of $\varphi \mid_{E(G)}$ has degree $n$, then*

$$[G : F^*(G)] \leq \mid I_{Out(E(G))}(\varphi \mid_{E(G)}) \mid < n(n+1) \leq \deg \varphi (\deg \varphi + 1).$$

*Proof.* We showed at the beginning of this section that

$$[G : F^*(G)] \leq \mid I_{Out(E(G))}(\varphi \mid_{E(G)}) \mid.$$

We apply Proposition 4. Unless $E(G)$ is of type $D_4$ or ${}^3D_4$,

$$\lambda a \log_m n \leq 2 \log_m n \leq 2 \log_2 n < n+1,$$

while, if $E(G)$ is of type $D_4$ or ${}^3D_4$, we have $m \geq 8$ and

$$\lambda a \log_m n \leq 6 \log_m n \leq 6 \log_8 n \leq 2 \log_2 n < n+1.$$

Finally, since $\varphi \mid_{E(G)}$ is homogeneous, $n$ divides $\deg \varphi$. $\square$

**Remark.** This Corollary is all we need for the proof of our main theorems. For large values of $n$, Proposition 4 gives a bound of order of magnitude $O(n \log n)$.

## 5. The proof of Theorem B

Throughout this section, we will assume the following hypothesis.

**Hypothesis I.** $G$ is a finite group such that $F^*(G) = Z(G).E^*(G)$, with an embedding $G \hookrightarrow GL(n,k)$ in which every proper normal subgroup of $G$ acts homogeneously.

Here, $E^*(G)$ is a central product of components and quasicomponents, the latter as defined by Proposition 2. The conditions are a consequence of primitivity, but we do not assume this; however, they still imply that $Z(G)$ is cyclic. Also, a component $E$ of $G$ will satisfy Hypothesis I (possibly for a smaller $n$) since every proper normal subgroup of $E$ will be characteristic in $G$, yet $E$ need not be primitive even when $G$ is[3].

Our next theorem appeals to the classification of finite simple groups for its exhaustive analysis.

**Theorem 6.** *Assume Hypothesis I with $E^*(G)$ quasisimple and irreducible. Then $[G : Z(G).E_\ell(G)]$ is bounded.*

*Suppose that $G = E(G) = E$, and let $\rho$ denote the representation of $E$ on the underlying vector space. If $[E : Z(E)].|I_{Out_\ell(E)}(\rho)| > (n+1)!$, then $E/Z(E) \notin \mathcal{L}_\ell$ and either*

(i)    $\ell$ *divides* $(n+2)$ *and* $E \cong A_{n+2}$, *or*

(ii)    $n \leq 12$ *and $E$ is one of the groups listed in* (a) – (c).

   (a)  *For almost all $\ell$ - not dividing $|Z(E)|$ and $E/Z(E) \notin \mathcal{L}_\ell$:*

---

[3]    For example, the alternating group $A_6$ has a 10-dimensional irreducible representation in every characteristic except 2 and 3 which is even monomial. This representation extends to primitive representations of $S_6$ and of $M_{10}$; the extension to a 10-dimensional representation of $PGL(2,9)$, however, remains monomial.

| $n$ | $E$ | $\mid I_{Out_{\ell}(E)}(\rho) \mid$ for most $\ell$ | $[E:Z(E)]$ | $\mid I_{Out_{\ell}(E)}(\rho) \mid$ for exceptional $\ell$ |
|---|---|---|---|---|
| 2 | $2 \cdot A_5$ | 1 | 60 | |
| 3 | $3 \cdot A_6$ | 1 | 360 | 2 when $\ell = 5$ |
| | $A_5$ | 1 | 60 | |
| | $L_3(2)$ | 1 | 168 | |
| 4 | $2 \cdot PSp_4(3)$ | 1 | 25920 | |
| | $2 \cdot L_3(2)$ | 1 | 168 | |
| | $2 \cdot A_6$ | 1 | 360 | 2 when $\ell = 5$ |
| | $2 \cdot A_7$ | 1 | 2520 | 2 when $\ell = 7$ |
| 5 | $PSp_4(3)$ | 1 | 25920 | |
| 6 | $6_1 \cdot U_4(3)$ | 2 | 3265920 | |
| | $U_3(3)$ | 2 | 6048 | |
| | $6 \cdot L_3(4)$ | 2 | 20160 | |
| | $PSp_4(3)$ | 2 | 25920 | |
| | $2 \cdot J_2$ | 1 | 604800 | 2 when $\ell = 5$ |
| 7 | $Sp_6(2)$ | 1 | 1451520 | |
| 8 | $2 \cdot O_8^+(2)$ | 2 | 174184400 | |
| | $2 \cdot Sp_6(2)$ | 1 | 1451520 | |
| 12 | $6 \cdot Suz$ | 1 | 448345497600 | |

(b) If $\ell$ divides $\mid Z(E) \mid$ in (a), the same results hold after factoring out $Z_\ell$, except for the cases $2 \cdot L_3(4)$ and $2 \cdot Suz$ when $\ell = 3$ that are included in (c).

(c) For particular values of $\ell$:

| $n$ | $\ell$ | $E$ | $\mid I_{Out_\ell(E)}(\rho) \mid$ | $[E:Z(E)]$ |
|---|---|---|---|---|
| 3 | 5 | $3 \cdot A_7$ | 1 | 2520 |
| 4 | 3 | $4_2 \cdot L_3(4)$ | 2 | 20160 |
| 5 | 3 | $M_{11}$ | 1 | 7920 |
| 6 | 3 | $2 \cdot L_3(4)$ | 4 | 20160 |
| | 2 | $3 \cdot M_{22}$ | 1 | 443520 |
| | 3 | $2 \cdot M_{12}$ | 1 | 95040 |
| 7 | 11 | $J_1$ | 1 | 175560 |
| 8 | 5 | $2 \cdot A_{10}$ | 1 | 1818400 |
| 9 | 2 | $3 \cdot J_3$ | 1 | 50232960 |
| 12 | 3 | $2 \cdot Suz$ | 2 | 448345497600 |

*Proof.* Put $E = E(G)$. If $E / Z(E) \in \mathcal{L}_\ell$, then $Z(G).E_\ell(G) = F^*(G)$ and the result holds by Corollary 5, so we may assume otherwise. Let $\rho$ denote the representation of $E$ on the underlying vector space. Then $G / E.Z(G)$ embeds into $I_{Out_\ell(E)}(\rho)$ so that we can establish boundedness by showing that only finitely many quasisimple groups remain that have faithful representations of degree $n$. So it is sufficient to determine all groups $G$ for which $G = E_{\ell'}(G) = E$ and $[E : Z(E)] \cdot | I_{Out_\ell(E)}(\rho) | > (n+1)!$

If $E$ is of Lie type, then $n \leq 8$ and only the groups listed in (a) can arise by Proposition 10 of [**3**]. All these groups are covered by the ATLAS [**5**] and the modular atlas [**8**], and $| I_{Out_\ell(E)}(\rho) |$ can be read off from the tables they contain[4]. If $E$ is sporadic, then we refer to [**7**] for the minimal degrees of projective representations in nonzero characteristic and simply check. In particular, $n \leq 12$. Leaving alternating groups aside for the moment, the only groups that can arise and whose representations do not arise by reduction $\mod \ell$ of a representation (possibly of a central extension of order $\ell$) in characteristic 0 and so do not already appear in Theorem 8 of [**3**] are given by the following table.

---

[4] The Introduction to the ATLAS explains how to use the tables. We note four particular features to explain our conclusions. Let $\alpha \in Aut(E / Z(E))$.

(i) $\alpha$ may permute irreducible representations of a covering group.

(ii) $\alpha$ may extend to an automorphism of $E$ but invert $Z(E)$.

(iii) If $\alpha$ has order 3 and $E / Z(E)$ has a central extension by $Z_2 \times Z_2$, then $\alpha$ may extend to this central extension and act nontrivially on the centre; then a cyclic extension does not admit any extension by $\alpha$.

(iv) If the irrationalities of a character occur just for one prime $p$, then for $\ell = p$ those characters become equal on reduction $\mod \ell$ and $| I_{Out_\ell(E)}(\rho) |$ may increase (the "exceptional $\ell$").

| $n$ | $[E:Z(E)] \cdot \mid I_{Out_\ell(E)}(\rho) \mid$ | $\ell$ | $G$ (for $\mid Z(G) \mid$ minimal) |
|---|---|---|---|
| 4 | 40320 | $\ell = 3$ | $4_2 \cdot (L_3(4).2_2)$ |
| 5 | 7920 | $\ell = 3$ | $M_{11}$ |
| 6 | 443520 | $\ell = 2$ | $3 \cdot M_{22}$ |
|   | 95040 | $\ell = 3$ | $2 \cdot M_{12}$ |
| 7 | 175560 | $\ell = 11$ | $J_1$ |
| 9 | 50232960 | $\ell = 2$ | $3 \cdot J_3$ |
| 12 | 896690995200 | $\ell = 3$ | $2 \cdot (Suz.2)$ |

We note from the final column of this table that, in all cases, there are groups for which $[G:Z(G)] = [E:Z(E)] \cdot \mid I_{Out_\ell(E)}(\rho) \mid$.

Suppose now that $E / Z(E) \cong A_m$. Dickson [**6**] showed that, for $m \geq 9$, the minimal degree faithful representation of $A_m$ in nonzero characteristic is $m - 1$ unless $\ell$ divides $m$, in which case it is $m - 2$, and that these representations extend to $S_m$. In this case, we get (i) immediately. If $m \leq 8$, when $\ell = 2$ we have the 4-dimensional representations of $A_7$ (that do not extend to $S_7$) – other alternating groups that might appear double as classical groups and lie in $\mathcal{L}_\ell$, while the outer automorphisms of $S_6$ do not stabilise the 5-dimensional representation of $A_6$ in any characteristic.

Suppose that $E(G) \cong 2 \cdot A_m$. Then we may assume that $\ell$ is odd. If $m > 7$, by a theorem of Wagner [**14**], the minimal degree $d$ of a faithful representation of $2 \cdot A_m$ is divisible by

$$2^{\left[\frac{m-s-1}{2}\right]}$$

where $s$ is the number of nonzero terms in the dyadic expansion of $m$. In particular, $s \leq \log_2(m+1)$ and this forces

$$d \geq \left(\frac{2^{m-3}}{m+1}\right)^{\frac{1}{2}} \geq 2^{\frac{m}{4}} \geq m$$

when $m \geq 16$, so that these cases do not arise; nor do they if $13 \leq m \leq 15$, when Wagner's result gives $d$ divisible by 16.

For $m \leq 12$ and every proper cover of $A_m$, all degrees are given explicitly in [8] (or the ATLAS for $\ell > m$), and we achieve the inequality $[E : Z(E)] \cdot |I_{Out_c(E)}(\rho)| > (n+1)!$ only when

$$E \cong \begin{cases} 2 \cdot A_5, \ n = 2, \ \ell \neq 2,5 \\ 2 \cdot A_6, \ n = 4, \ \ell \neq 2,3 \\ 3 \cdot A_6, \ n = 3, \ \ell \neq 2,3 \\ 2 \cdot A_7, \ n = 4, \ \ell \neq 2 \\ 3 \cdot A_7, \ n = 3, \ \ell = 5 \\ 2 \cdot A_{10}, \ n = 8, \ \ell = 5. \end{cases}$$

where we have again excluded cases arising from isomorphisms with classical groups. □

We will need the following counterpart for quasicomponents.

**Theorem 7.** *Assume Hypothesis I with $n = p^m$ and $E^*(G)$ a quasicomponent of order $p^{2m+1}$. Suppose that $[G : Z(G)]$ is maximal. Then $[G : Z(G)]$ is greater than the maximal index in Theorem 6 for the pair $(n, k)$ only in the following cases, where the group $G$ stated satisfies $|Z(G)| = p$:*

| $n$ | $p$ | $\ell$ | $[G : Z(G)]$ | $G$ |
|-----|-----|--------|--------------|-----|
| 2 | 2 | 5 | 24 | $GL(2,3)$ |
| 3 | 3 | 2 | 216 | $3^{1+2}.SL_2(3)$ |
| 5 | 5 | 2 | 3000 | $5^{1+2}.SL_2(5)$ |
| 9 | 3 | $\ell \neq 3$ | 4199040 | $3^{1+4}.Sp_4(3)$ |

*Proof.* By Proposition 2(iv), $G / F^*(G)$ is isomorphic to a subgroup of $Sp_{2m}(p)$. By Theorem 6 of [3] and the subsequent remark, we obtain a bound for $[G : Z(G)]$ of $p^{2m}|Sp_{2m}(p)|$, and that $p^n = 2, 3, 4, 5, 8$ or $9$. We must have $p \neq \ell$. If $p^m = 2$, then $G / Z(G) \cong S_4$, and this will arise when $\ell = 5$. Careful comparison of orders with those of the groups occurring in Theorem 6 eliminates the cases $p^m = 4$ or $8$, leaving the possibilities $p^m = 3, 5$, or $9$, which can occur. These give bounds of 216, 3000 and 4199040 respectively, for $\ell = 2$ in the first two cases and $\ell \neq 3$ in the last. By Proposition 7 of [**3**], these are achieved by the semidirect product of an extraspecial group of exponent $p$ and a symplectic group, and the representations are the $\mod \ell$ reductions of characteristic 0 representations. (Note that $Sp_2(p) \cong SL_2(p)$.) $\square$

We now embark on the proof of Theorem B. First we note the following.

**Proposition 8.** *Except for the case $n = 9$ and $\ell = 5$, the bounds claimed in Theorem B(ii) are the maximum of the values for $[E : Z(E)] \cdot | I_{Out_\ell(E)}(\rho) |$ given by Theorem 6 or for $[G : Z(G)]$ given in Theorem 7.*

Let $G$ be a primitive subgroup of $GL(n,k)$. Let $V$ be the associated $kG$-module and let $U$ be a simple $kE^*(G)$-submodule. Let $E_1, \ldots, E_s$ be the components and quasicomponents of $G$ (no longer distinguished). As in the proof of Theorem 3, we may decompose $U$ as a tensor product

$$U = U_1 \otimes \cdots \otimes U_s$$

where $U_i$ is a simple $kE_i$-module affording a representation $\rho_i$ of $E_i$ of degree $n_i$ with $\ker(\rho_i) \subseteq Z(E_i)$ for each $i$. We call $n_1, \ldots, n_s$ the *subdegrees* corresponding to the (quasi)components of $G$. Then $n_1 \cdots n_s \mid n$ and, by Theorem 3 and Corollary 5,

$$(*) \qquad [G : Z(G).E_\ell(G)] \leq b_1 \cdots b_s \prod_{j=1}^{t} l_j !$$

where $b_1, \ldots, b_s$ are the *bounds for the contributions* of the (quasi)components defined by

$$b_i = \begin{cases} p_i^{2m_i} \mid Sp_{2m_i}(p_i) \mid & \text{if } E_i \text{ is a quasicomponent with } [E_i : Z(E_i)] = p_i^{2m_i} \\ [E_i : Z(E_i)] \cdot | I_{Out_\ell(E_i)}(\rho_i) | & \text{if } E_i \notin \mathcal{L}_\ell \\ n_i(n_i + 1) & \text{if } E_i \in \mathcal{L}_\ell, \end{cases}$$

and $l_1, \ldots, l_t$ are the sizes of the distinct isomorphism classes of components. For given $n$, we can apply Theorems 6 and 7 to each component or quasicomponent in turn to show that the values $b_i$ are bounded. So the index $[G : Z(G).E_\ell(G)]$ is bounded.

Our goal now will be to show that in general this index is maximal when there is exactly one component or quasicomponent. There is just the one exception that arose in Proposition 8.

**Lemma 9.** *Suppose that $\ell = 5$ and let $G$ be a primitive subgroup of $GL(n, k)$ with $[G : Z(G).E_\ell(G)]$ maximal. Then no component or quasicomponent of $G$ has degree 9. If $n = 9$, then $G / Z(G) \cong A_7 wr S_2$.*

*Proof.* By Theorems 6 and 7, the maximal contribution of a component or quasicomponent of degree 9 to the index $[G : Z(G).E_\ell(G)]$ is 4199040. Since $|A_7|^2 = 6350400$, we can increase the index by replacing such a component by the central product $(3A_7) \circ (3A_7)$. So there is no (quasi)component of degree 9 when $[G : Z(G).E_\ell(G)]$ is maximal.

If $n = 9$, we must then have $E(G) \cong (3A_7) \circ (3A_7)$ with the two components interchanged by an element of $G \setminus F^*(G)$, and the bound claimed in Theorem B follows. $\square$

In the light of Proposition 8 and Lemma 9, we define, for any positive integer $r$ and each $\ell$, a constant $N_{r,\ell}$ to be the bound for degree $r$ and characteristic $\ell$ where it was claimed in Theorem B(ii) and by

$$N_{r,\ell} = \begin{cases} (r+2)! & \text{if } r \geq 8, \ r \neq 12, \text{ and } \ell \text{ divides } r+2 \\ (r+1)! & \text{otherwise.} \end{cases}$$

We note that, except when $\ell$ divides $(r+1)$ and for $N_{9,5}$, the value $N_{r,\ell}$ can be achieved as the bound for the contribution of some component or quasicomponent in a primitive group. Furthermore, that bound is achieved, either by a symmetric group or by a group listed explicitly in part (ii) of Theorem B. In the bound (*) for the index $[G : Z(G).E_\ell(G)]$, we may replace each $b_i$ by the $N_{r,\ell}$ corresponding to the appropriate subdegree to get

$$[G : Z(G).E_\ell(G)] \leq \prod_{j=1}^{t} (N_{r_j,\ell})^{l_j} l_j!$$

where now $r_1, \ldots, r_t$ are the distinct subdegrees and $l_1, \ldots, l_t$ are the numbers of (quasi)components of each subdegree so that $r_1^{l_1} \cdots r_t^{l_t}$ divides $n$.

Suppose now that $n$ and $\ell$ are fixed. We can work under the following hypothesis, noting that the subsequent arguments do not require us to distinguish the cases of contributions that arise from quasicomponents by insisting then that $l_j = 1$. If we can show that maximality forces $t = 1$ and $l_1 = 1$, then we can investigate the consequences.

**Hypothesis II.** The distinct subdegrees $r_1, \ldots, r_t$ and their multiplicities $l_1, \ldots, l_t$ are chosen so that $r_1^{l_1} \cdots r_t^{l_t}$ divides $n$ and

$$\prod_{j=1}^{t} (N_{r_j,\ell})^{l_j} l_j!$$

is maximal.

As in [**3**], we will refer to components when analysing this situation, even though no group may exist.

**Lemma 10.** *The inequality* $(N_{r,\ell})^{l_r} \cdot (l_r)! < (r^L + 1)!$ *holds when*

(i)  $l_r > 1$ *for* $r \geq 3$*, except when* $\ell = 5$*,*

(ii)  $l_r > 1$ *for* $r \geq 4$ *or* $l_r > 2$ *for* $r = 3$*, when* $\ell = 5$*,*

(iii)  $l_r > 3$ *for* $r = 2$ *if* $\ell \neq 2, 5$*,*

(iv)  $l_r > 2$ *for* $r = 2$ *if* $\ell = 5$*, or*

(v)  $l_r > 1$ *for* $r = 2$ *if* $\ell = 2$*.*

*Proof.* If $(N_{r,\ell})^L \cdot (L)! < (r^L + 1)!$ for some $L > 1$, then

$$(N_{r,\ell})^{L+1} \cdot (L+1)! < ((N_{r,\ell})^L \cdot L!)^2 < ((r^L+1)!)^2 < (r^{L+1}+1)!$$

so that the inequalities will hold once the base cases are established for each pair $r, \ell$. If $r \geq 3$, then $2 \cdot (r+2)! < (r^2+1)!$, establishing the desired inequality in the "generic" situation for all $\ell$; for the exceptional values of $N_{r,\ell}$ when $3 \leq r \leq 12$, only the case $r = 3$, $\ell = 5$, $l_3 = 2$, $N_{3,5} = 2520$ fails. If $r = 2$, only the values $l_r = 2$ and $l_r = 3$ provide exceptions. $\square$

**Corollary 11.** *Assume Hypothesis II. Then we may suppose that $l_j = 1$ for all $j$.*

*Proof.* For all $r$, $L$ and $\ell$, $(r^L+1)! \leq N_{r^L,\ell}$. Hence, when Lemma 10 applies, we can replace any $L$ components of degree $r$ by a single component of degree $r^L$, increasing the multiplicity of such components by 1. But this would contradict the maximality. The only obstructions to this process occur when $(r, L) = (3, 2)$ if $\ell = 5$, and for the cases $r = 2$, $L \leq 3$. In the latter cases, we can check that $2 \cdot (N_{2,\ell})^2 < N_{4,\ell}$ and $6 \cdot (N_{2,\ell})^3 < N_{8,\ell}$ for all $\ell$ so that we can carry out the same replacement. In the remaining case, we note that $2 \cdot (N_{3,5})^2 = N_{9,5}$, and we may make that substitution. $\square$

**Lemma 12.** *Assume Hypothesis II and that $l_j = 1$ for all $j$. Then $t = 1$.*

*Proof.* We first show that $N_{p,\ell} N_{q,\ell} < (pq+1)!$ whenever $p < q$ and $pq > 12$, noting that $(pq+1)! \leq N_{pq,\ell}$ for all $\ell$. This can be seen by direct calculation if $q \leq 12$ and, when $p \leq 12$ and $q \geq 13$, from the inequalities

$$\frac{(pq+1)!}{N_{q,\ell}} \geq \frac{(2q+1)!}{N_{q,\ell}} \geq (2q+1)\cdots(q+3) \geq (16)^{12} = 2^{48} > 2^8 \cdot 10^{12} > N_{p,\ell}.$$

If $13 \leq p < q$, then a similar argument yields the inequality

$$(p+2)! \cdot (q+2)! < (pq+1)!$$

For $pq \leq 12$, we simply check that $N_{p,\ell} N_{q,\ell} < N_{pq,\ell}$ for $p \neq q$ to make the analogous substitution. $\square$

Thus we have shown that, if $G$ is a primitive finite subgroup of $GL(n,k)$, then

$$[G : Z(G).E_\ell(G)] \leq N_{r,\ell}$$

for some $r$ dividing $n$. But now, except when $\ell$ divides $n+1$, the existence of a primitive group of order $N_{n,\ell}$ – either a symmetric group or a group listed in the conclusion of Theorem B – establishes the theorem. In the case that $\ell$ does divide $n+1$, we have $N_{n,\ell} = (n+1)!$, and it is easy to verify that $N_{n,\ell} > N_{r,\ell}$ for every proper divisor $r$ of $n$, so that $(n+1)!$ *is* a bound, though not one that is necessarily achieved.

## 6. Pairs and a replacement theorem

We will use our bounds for primitive groups to establish Theorem A. We need to generalise the methods of [**4**] to take account of the subgroup $E_\ell(G)$; we do so in slightly greater generality than needed, and in this section prove a technical result that will allow us to isolate any family of components.

Let $G$ be a finite subgroup of $GL(n,k)$. By Proposition 1, when proving Theorem A, we may suppose that $O_\ell(G) = 1$ and that $G$ acts completely reducibly on the underlying vector space. Thus we shall suppose for the remainder of this paper that $G$ satisfies the following hypothesis.

**Hypothesis III.** $G$ acts faithfully on a vector space $V$ of dimension $n$, and there is a decomposition

$$V = V_1 \oplus \cdots \oplus V_r$$

of $V$ whose summands are permuted by $G$ and such that, if $H_i = Stab_G(V_i)$, then $H_i$ acts primitively (but not necessarily faithfully) on $V_i$ for $1 \leq i \leq r$.

**Definition 13.** A *primitive decomposition pair* $(G,V)$ consists of a group $G$ and a vector space $V$ together with a decomposition of $V$ for which Hypothesis III holds. We will write $(G,k^n)$ or $(G,V;n,k)$ if we wish to emphasise the dimension $n$ of $V$ or the field $k$ (in particular, its characteristic $\ell$), and talk just of the *pair* $(G,V)$. A pair will be *primitive* if $G$ acts primitively on $V$ (i.e., $r = 1$).

We will need to study certain distinguished pairs.

**Definition 14.** (i) Let $(G_1,V_1),\ldots,(G_t,V_t)$ be a collection of pairs. Then their *sum* is the pair $(\tilde{G},\tilde{V})$ where $\tilde{G} = G_1 \times \cdots \times G_t$, $\tilde{V} = V_1 \oplus \cdots \oplus V_t$ as a vector space, made into a $k\tilde{G}$-module via

$$(g_1,\ldots,g_t)(v_1 + \cdots + v_t) = g_1 v_1 + \cdots + g_t v_t \qquad (g_i \in G_i, v_i \in V_i)$$

and we write $(\tilde{G},\tilde{V}) = \sum_{i=1}^{t}(G_i,V_i)$.

(ii) If $(G_i,V_i) = (G,V)$ for all $i$, put $(\tilde{G},\tilde{V}) = (G,V)^t = \sum_{i=1}^{t}(G_i,V_i)$.

(iii) Given (ii), we can extend $\tilde{V}$ to a $k(G\,wr\,S_t)$-module $\hat{V}$ via

$$\sigma^{-1}(g_1,\ldots,g_t)\sigma = (g_{\sigma^{-1}1},\ldots,g_{\sigma^{-1}t}) \quad \text{and} \quad \sigma(v_1 + \cdots + v_t) = v_{\sigma 1} + \cdots + v_{\sigma t},$$

and define the *wreath product* $(G,V)wr\,S_t = (G\,wr\,S_t,\hat{V})$.

(iv) A pair $(G,V)$ is said to be *saturated* if

$$(G,V) = \sum_{i=1}^{s}(G_i,V_i)wr\,S_{t_i}$$

where each pair $(G_i,V_i)$ is primitive.

**Remarks.** (i) $(\tilde{G},\tilde{V})$ as defined in (i) *is* a primitive decomposition pair since it inherits a primitive decomposition from each summand $V_i$.

(ii)   The wreath product construction in (iii) does give rise to a primitive decomposition pair since the primitive decomposition of (ii) is preserved by the action of the symmetric group $S_t$.

Let $(G, V)$ be a primitive decomposition pair, and let

$$H = H(G) = \bigcap_{i=i}^{r} H_i \, .$$

Then $H$ is the kernel of the permutation action of $G$ on the set of subspaces $\{V_i\}$ so that, in particular, $H \lhd H_i$ for each $i$. In the notation of Hypothesis III, let $K_i$ be the kernel of the action of $H_i$ on $V_i$ and put $P_i = H_i / K_i$. Then we may regard $V_i$ as a faithful $kP_i$-module. In the language of pairs, we may restate Theorem 2 (the replacement theorem) of [**4**] in the following way.

**Theorem 15.**   *Let $(G, V)$ be a primitive decomposition pair and assume the notation above. Then, for a suitable ordering of the summands $\{V_i\}$, there is a saturated pair*

$$(\hat{G}, \hat{V}) = \sum_{i=1}^{s} (P_i, V_i) wr S_{t_i}$$

*where $\dim \hat{V} = \dim V$, and a natural embedding $H \hookrightarrow \hat{H} = P_1 \times \cdots \times P_r$ for which $Z(H) = H \cap Z(\hat{H})$. Furthermore, $[G : Z(H)] \leq [\hat{G} : Z(\hat{H})]$.*

It is clear from the definitions that $\hat{H}$ is the kernel of the permutation action of $\hat{G}$ on the summands of $\hat{V}$ so that the inequality $[G : Z(H)] \leq [\hat{G} : Z(\hat{H})]$ is actually a consequence of the property $Z(H) = H \cap Z(\hat{H})$, which was contained in the proof. We now need to modify this.

**Definition 16.**   Let $\mathfrak{L}$ be any collection of finite simple groups. For an arbitrary finite group $X$, put

$$E_{\mathfrak{L}}(X) = \prod_{\substack{E \vartriangleleft E(X) \\ E/Z(E) \in \mathfrak{L}}} E \, ,$$

namely, the product of those components of $X$ whose simple quotients lie in $\mathfrak{L}$.

**Remark.** In the statement of Theorem A, we put $E_\ell(G) = E_{\mathcal{L}_\ell}(G)$.

The following technical result may be of independent interest.

**Proposition 17.** *Let* $X = X_1 \times \cdots \times X_n$ *be a direct product of nonabelian finite simple groups and let* $Y$ *be a subgroup whose projections onto each simple direct factor of* $X$ *are either surjective or trivial. Then* $Y = E(Y)$. *Furthermore, for any family* $\mathfrak{L}$ *of finite simple groups,* $E_{\mathfrak{L}}(Y) = Y \cap E_{\mathfrak{L}}(X)$.

*Proof.* Let $\pi_i : X \to X_i$ be the natural projection for each $i$. If $E \vartriangleleft\vartriangleleft Y$, then either $\pi_i(E) = X_i$ or $\pi_i(E) = 1$. Thus $F(Y) = 1$ and $F^*(Y) = E(Y)$; in particular, the components of $Y$ are also simple.

If $Y \neq E(Y)$, pick $y \in Y \setminus E(Y)$. Then there is a component $E$ of $Y$ for which either $E^y \neq E$ or $y$ acts on $E$ as an outer automorphism. In either case, $E(\langle y, E \rangle)$ is the unique minimal normal subgroup of $\langle y, E \rangle$ and contains $E$. Now $\pi_i(E) \neq 1$ for some $i$, and $\pi_i\big|_{\langle y, E \rangle}$ is injective since $E(\langle y, E \rangle) \cap \ker \pi_i = 1$; hence

$$\pi_i(\langle y, E \rangle) \supset \pi_i(E) = X_i \, ,$$

which is impossible. So $Y = E(Y)$.

Now $E \subseteq E_{\mathfrak{L}}(X)$ for any component $E$ of $E_{\mathfrak{L}}(Y)$ so that $E_{\mathfrak{L}}(Y) = Y \cap E_{\mathfrak{L}}(X)$. $\square$

Our new replacement theorem is the following.

**Theorem 18.** *The final inequality of Theorem 15 may be replaced by*

$$[G : Z(H)E_{\mathcal{L}}(H)] \leq [\hat{G} : Z(\hat{H})E_{\mathcal{L}}(\hat{H})]$$

*for any collection $\mathcal{L}$ of finite simple groups.*

*Proof.* Since $[G : H] \leq [\hat{G} : \hat{H}]$, we need only to establish the inequality

$$[\hat{H} : Z(\hat{H})E_{\mathcal{L}}(\hat{H})] \geq [H : Z(H)E_{\mathcal{L}}(H)],$$

and we prove this entirely within $\hat{H}$. Thus we may ignore the wreath product construction and write $\hat{H} = P_1 \times \cdots \times P_r$, and identify $H$ with its embedding, which came via the homomorphisms $H \hookrightarrow H_i \rightarrow P_i$.

For any group $X$,

$$E_{\mathcal{L}}(X / Z(X)) = E_{\mathcal{L}}(X).Z(X) / Z(X);$$

Applying this both to $\hat{H}$ and to $H$, since $Z(H) = H \cap Z(\hat{H})$ we may assume without loss that $Z(\hat{H}) = Z(H) = 1$, and then need only show that

$$E_{\mathcal{L}}(H) = H \cap E_{\mathcal{L}}(\hat{H}).$$

Since $H \lhd H_i$, we have $\pi_i(H) \lhd P_i$ for each $i$ where $\pi_i : \hat{H} \rightarrow P_i$ is the natural projection. Let $E$ be a component of $H$. Then $\pi_i(E) \lhd\lhd P_i$ for each $i$, and it follows that $\pi_i(E)$ is either trivial or a component of $P_i$. In particular, $E(H) \subseteq E(\hat{H})$.

Since $Z(\hat{H}) = 1$, $E(\hat{H})$ is a direct product of (uniquely determined) simple subgroups. Now $\pi_i(H \cap E(\hat{H})) \subseteq \pi_i(E(\hat{H})) = E(P_i) \lhd E(\hat{H})$; since $H \cap E(\hat{H}) \lhd H$, in addition $\pi_i(H \cap E(\hat{H})) \lhd\lhd P_i \lhd \hat{H}$ so that $\pi_i(H \cap E(\hat{H})) \lhd\lhd E(\hat{H})$ and further $\pi_i(H \cap E(\hat{H})) \lhd E(\hat{H})$. Thus, if $\hat{E}$ is a component of $\hat{H}$, then $H \cap E(\hat{H})$ projects into $\hat{E}$ either surjectively or trivially.

The hypothesis of Proposition 17 is now satisfied with $X = E(\hat{H})$ and $Y = H \cap E(\hat{H})$. Hence, as $H \cap E(\hat{H}) \lhd H$ so that every component of $H \cap E(\hat{H})$ is also a component of $H$,

$$H \cap E(\hat{H}) = E(H \cap E(\hat{H})) \subseteq E(H) \subseteq E(\hat{H})$$

from which, taking intersections with $H$, we first deduce that $E(H) = H \cap E(\hat{H})$. Then, applying the latter part of Proposition 17 (with $Y = E(H)$ now), we obtain

$$E_{\mathfrak{L}}(H) = E_{\mathfrak{L}}(E(H)) = E(H) \cap E_{\mathfrak{L}}(\hat{H}) = H \cap E_{\mathfrak{L}}(\hat{H}),$$

as required. □

## 7. The proof of Theorem A and the extension to Theorem A$'$

Let $G$ be a finite subgroup of $GL(n,k)$. Then, by Theorem 18, there is a saturated pair

$$(\hat{G}, \hat{V}; n, k) = \sum_{i=1}^{s} (P_i, V_i) wr S_{t_i}$$

such that $[G : Z(H)E_\ell(H)] \leq [\hat{G} : Z(\hat{H})E_\ell(\hat{H})]$. Now

$$[G : Z(H).E_\ell(G)] \leq [G : Z(H).E_\ell(H)] \leq [\hat{G} : Z(\hat{H}).E_\ell(\hat{H})]$$
$$= \prod_{i=1}^{s} ([P_i : Z(P_i).E_\ell(P_i)])^{t_i} t_i!$$

By Theorem B, each term $[P_i : Z(P_i).E_\ell(P_i)]$ in this product is bounded in terms of $\dim \hat{V}_i$ and hence, as

$$n = \sum_{i=1}^{s} t_i \dim \hat{V}_i,$$

$[G : Z(H).E_\ell(G)]$ is bounded.

In order to establish an actual bound, it is sufficient to consider saturated pairs only. So, for the remainder of the paper, for each $k$ and each degree $n$ we pick a saturated pair

$$(G, V; n, k) = \sum_{i=1}^{s} (P_i, V_i) wr S_{t_i}$$

for which, with $H = H(G)$ as defined prior to Theorem 15,

(I)    $[G : Z(H)E_\ell(G)]$ is maximal.

Then, if we regard $V$ as a $kH$-module and write

$$(H,V) = \sum_{i=1}^{r} (P_i, V_i),$$

we will call each subgroup $P_i$ (*resp.* pair $(P_i, V_i)$) a (*primitive*) *constituent* of $G$ (*resp.* $(G,V)$) and $\dim V_i$ the corresponding *subdegree*.

We note that $E_\ell(G) = E_\ell(H) = E_\ell(P_1) \times \cdots \times E_\ell(P_r)$ for a saturated pair $(G,V)$ and also that $Z(G) = Z(H) = Z(P_1) \times \cdots \times Z(P_r)$; thus

$$[G : Z(H).E_\ell(G)] = \prod_{i=1}^{s} ([P_i : Z(P_i).E_\ell(P_i)])^{t_i} t_i !$$

This formula, with the maximality of $[G : Z(H)E_\ell(G)]$, now forces the following.

**Lemma 19.** *If $(P_i, V_i)$ and $(P_j, V_j)$ are primitive constituents of the same subdegree, then $P_i \cong P_j$ and $V_i \cong V_j$ as $kP_i$-modules, and exactly one of $P_1, \ldots, P_r$ has any given subdegree.*

*Proof.* We first require $[P_i : Z(P_i).E_\ell(P_i)]$ maximal, given $\dim V_i$, for each constituent in order to maximise the index $[G : Z(H)E_\ell(G)]$. But now maximality is achieved only if, when $(P_i, V_i)$ and $(P_j, V_j)$ have the same subdegree, they appear within the same wreath product, forcing the conclusion. $\square$

The maximality of $[P_i : Z(P_i).E_\ell(P_i)]$ without further restriction means that we can choose $P_i$ whenever we know a group that attains the bound. It is convenient also to assume

(II) Each primitive constituent $(P_i, V_i)$ is chosen such that, given $\dim V_i$, first the index $[P_i : Z(P_i).E_\ell(P_i)]$ is maximal, and then $|Z(P_i)|$ is minimal.

In general, the choice can be made from the groups listed in Theorem B (including symmetric groups), but a problem arises when $\ell$ divides $n_i + 1$ and we have only the bound $(n_i + 1)!$ for the index; in this case the symmetric group $S_{n_i+1}$ has only a reducible representation of degree $n_i$, but does have a primitive representation of degree $n_i - 1$. In this situation, even if there is a potential primitive constituent of degree $n_i$ and order $(n_i + 1)!$, we replace every such constituent with the symmetric group $S_{n_i+1}$ as a primitive group of degree $n_i - 1$ together with a primitive constituent of degree 1. With this choice, and making a similar replacement when there are larger primitive groups of smaller degree, we now show that primitive constituents of our saturated pair $(G, V)$ for which (I) holds can be chosen in the following way.

**Lemma 20.** *A primitive constituent $P$ of subdegree $m$ may be taken to be one of the following:*

(i)  $P = 1$ *when* $m = 1$;

(ii)  *if* $2 \le m \le 9$ *or* $m = 12$, *one of the groups listed in Theorem B(ii) except that*

  (a)  *there is no primitive constituent of subdegree 2 if* $\ell = 2$,

  (b)  *there is no primitive constituent of subdegree 5 or 9 unless* $\ell = 2$,

  (c)  *there is no primitive constituent of subdegree 7 unless* $\ell = 3$, *and*

  (d)  *there is no primitive constituent of subdegree 8 if* $\ell = 2$;

(iii)  *if* $m \ge 10, m \ne 12$, *then there is no primitive constituent if* $m \equiv -1 \pmod{\ell}$ *and*

$$P \cong \begin{cases} S_{m+2} & \text{if } m \equiv -2 \pmod{\ell} \\ S_{m+1} & \text{otherwise.} \end{cases}$$

*Proof.* Part (i) follows from the minimality of $|Z(P)|$. Now all the remaining claims except for (ii)(a) follow from Theorem B and a replacement argument similar to that already described for (iii) above. In the outstanding case, we would first take $P \cong SL(2, 2^a)$ but then replace the pair $(P, k^2)$ by $(1, k) wr S_2$, contradicting the maximality of $[G : Z(H) E_\ell(G)]$.  $\square$

As an immediate consequence of this explicit list of possible primitive constituents, and also by Lemma 19, we will change our notation and let $P_{m,\ell}$ stand for a primitive constituent of subdegree $m$ when $k$ has characteristic $\ell$ and let $t_{m,\ell}$ denote its multiplicity; when the characteristic $\ell$ is explicitly chosen, then we will write just $P_m$ for $P_{m,\ell}$, and $t_m$ (or just $t$) for $t_{m,\ell}$. We will also let $P_m$ denote the pair $(P_m, k^m)$ when there is no risk of confusion.

In addition, we have the following corollary.

**Corollary 21.** $E_\ell(G) = 1$. In particular, $[G : Z(H)E_\ell(G)] = [G : Z(H)]$.

From this point, having eliminated $E_\ell(G)$ from consideration in order to determine bounds, the argument can essentially follow the same route as in the characteristic zero case [**4**]. We therefore only indicate the significant differences and sketch the main steps since the details are purely numerical.

If $\ell \geq 7$, apart from the exceptional case in Lemma 20(iii) when $m \equiv -2 \pmod{\ell}$, all primitive constituents are precisely the same as those arising in characteristic zero. For $\ell \leq 5$, the maximal primitive constituents in small degree may be quite different. In any case, our goal is to show that the maximality of the index $[G : Z(H)]$ can be achieved for a saturated pair for large enough $n$ only when there is a single nontrivial primitive constituent, and then appeal (again) to Lemma 20 to show that it is a symmetric group $S_{n+1}$ or $S_{n+2}$. This reduction is carried out by a "replacement" process; we eliminate potential summands of our saturated pair by showing that there could then be different summands of the same total degree that would give a greater contribution to the index $[G : Z(H)]$, contrary to the maximality choice in (I).

First, refinements of Lemma 4 of [**4**] give, for $m, m' \geq 8$, the inequalities $((m+2)!)^t t! < (tm+1)!$ and $(m+2)!(m'+2)! < (m+m'+1)!$ so that the replacements of a direct factor $P_{m,\ell} wr S_t$ by a primitive pair $(P_{mt,\ell}, k^{mt})$ and any sum of pairs $(P_{m,\ell}, k^m) + (P_{m',\ell}, k^{m'})$ by the pair $(P_{m+m',\ell}, k^{m+m'})$ whenever $P_{m,\ell}$ and $P_{m',\ell}$ are symmetric groups give the following.

**Lemma 22.** *For any $\ell$, at most one primitive constituent is a symmetric group.*

The possibility that $P_m \cong S_{m+2}$ means that the arguments of Lemma 9 of [**4**] need to be modified; in particular, symmetric groups of smaller degrees may occur as primitive constituents here than occur in characteristic zero.

**Lemma 23.** *If $\ell \neq 2$ or $5$, then any primitive constituent that is a symmetric group has subdegree at least $53$.*

*Proof.* In smaller degree, we could replace $P_m$ by a wreath product $(2 \cdot A_5)wrS_t$ for some $t$. (In fact, this lemma is best possible for $\ell = 11$.) $\square$

**Lemma 24.** (i) *At most one primitive constituent has subdegree $1$ unless $\ell = 2$, in which case there may be two.*

(ii) *If any primitive constituent is a symmetric group with subdegree m, then there is none of subdegree $1$ unless $m \equiv -2 \pmod \ell$, and then at most one if $\ell = 2$.*

*Proof.* Replacing $(1,k)wrS_t$ by $(P_t, k^t)$ yields (i), while replacing $(P_{m,\ell}, k^m) + (1,k)wrS_2$ by $(P_{m+2,\ell}, k^{m+2})$ and $(P_{m,\ell}, k^m) + (1,k)$ by $(P_{m+1,\ell}, k^{m+1})$ yields (ii). $\square$

As we noted above, if $\ell \geq 7$, every primitive constituent is the same as in characteristic zero. We can now appeal to the calculations in [**4**]. Minor refinements to the numerical arguments of Lemma 9 of [**4**] now give our generic bounds, and that they fail for $n = 70$.

**Proposition 25.** *If $\ell \geq 7$, then the bounds of Theorems A and $A'$ hold.*

We now complete the proof of Theorem A in this case; the same will hold for smaller primes when the bound is achieved by a symmetric group.

**Proposition 26.** Let $G$ be a finite subgroup of $GL(n, k)$ for which the bound of Theorem A is achieved. Then the structure of $G$ is as claimed.

*Proof.* The saturated pair $(\hat{G}, \hat{V})$ constructed from $G$ via Proposition 1 and Theorem 15 is already maximal. There is one primitive constituent that is an appropriate symmetric group, together with just one trivial constituent when $n \equiv -1 \pmod{\ell}$. In either case, the construction of Theorem 15 yields an injection of $\bar{G}$ into $\hat{G}$, and maximality gives $E(\bar{G})$ simple and the isomorphism $\bar{G} / Z(\bar{G}) \cong \hat{G} / Z(\hat{G})$, where $\bar{G} = G / O_\ell(G)$.

In the nonexceptional cases, $G$ acts irreducibly and hence $O_\ell(G) = 1$. Otherwise the underlying vector space $V$ has two composition factors as a $kG$-module, one of which is trivial, and $G$ stabilises a flag $0 \subset V' \subset V$ where one factor has dimension 1 and $G$ acts irreducibly on the other; hence $O_\ell(G)$ is elementary abelian. $\square$

The situation when $\ell = 3$ is very similar, and we need only establish the bounds.

**Lemma 27.** If $\ell = 3$, there is no primitive constituent of subdegree $m$ for $4 < m \leq 12$.

*Proof.* If there were, and applying Lemma 20, successive replacements of $P_{12}$ by $P_4 wr S_3$, $P_{11}$ by $P_{10} + P_1$, $P_{10}$ by $P_4 wr S_2 + P_2$, $P_8$ by $P_4 wr S_2$, $P_7$ by $P_4 + P_3$, and $P_6$ by $P_4 + P_2$ would yield a contradiction. $\square$

**Proposition 28.** Theorem $A'$ holds for $\ell = 3$.

*Proof.* If $G = (2 \cdot A_5)wrS_{34}$, then $[G : Z(H)] > 69!$ so we may suppose that $n \geq 69$. However, with $\ell = 3$ we can replace $(2 \cdot A_5)wrS_t$ by a symmetric group whenever $t \geq 35$ so it suffices to eliminate primitive constituents of subdegree 4; then the same arguments as for characteristic zero (but with just subdegrees 2 and 3 to consider) apply. So suppose that there are $t$ such constituents. We can replace $P_4wrS_t$ by $(2 \cdot A_5)wrS_{2t}$ if $t \geq 8$; so we have $t \leq 7$. Similarly, there is at most one primitive constituent of subdegree 3.

If any primitive constituent is a symmetric group $S_m$, then Lemma 23 can be modified to show that $m \geq 57$ since 3 does not divide 55. Now we can replace $S_m + P_4wrS_t$ by a single symmetric group $S_{m+4t-1}$ (perhaps with an extra trivial constituent).

If not, then if $n \geq 70$ there must be at least 19 primitive constituents of subdegree 2, and then the replacement of $(2 \cdot A_5)wrS_r + P_twrS_t$ by $(2 \cdot A_5)wrS_{r+2t}$ yields a contradiction. $\square$

The cases where $\ell = 2$ and $\ell = 5$ are somewhat different and we treat them in turn.

## 8. The case $\ell = 2$

By Lemma 20, there are no primitive constituents of subdegrees 2, 7 or 8. Now the successive replacements of $P_{13}$ by $P_{12} + P_1$, $P_{12}$ by $P_6wrS_2$, $P_{11}$ by $P_{10} + P_1$, $P_{10}$ by $P_6 + P_4$, and of $P_9$ by $P_6 + P_3$ eliminate subdegrees $m$ for $6 < m \leq 13$.

Next, the replacement for $3 \leq m \leq 6$ of $P_mwrS_t$ by either $(S_{tm+1}, k^{tm})$ or $(S_{tm+2}, k^{tm})$ reduces the number of primitive constituents of subdegrees 3, 4, 5 or 6 to at most 5, 3, 1 or 5, respectively, and then the further replacements of $P_3wrS_t$ by $P_6$, $P_6 + P_3$, $P_6wrS_2$ or $P_6wrS_2 + P_3$ for multiplicities 2, 3, 4 or 5 and of $P_4wrS_t$ by $P_6wrS_2$ or $P_6 + S_2$ for multiplicities 2 or 3 reduce the multiplicities of $P_3$ and $P_4$ to at most one each.

**Lemma 29.** *At most one primitive constituent has subdegree $d$ with $3 \leq d \leq 5$. If there is any, there can be no primitive constituent of degree 1.*

*Proof.* We could replace $P_3 + P_4 + P_5$ by $P_6 wr S_2$, any $P_m + P_{m'}$ by a sum of $P_6$ and trivial constituents, and any $P_m + P_1$ by $P_{m+1}$. □

**Lemma 30.** *If any primitive constituent is a symmetric group, then there is no other nontrivial primitive constituent.*

*Proof.* Such a constituent $P_m$ is a symmetric group $S_q$ for some $q \geq 15$, by Lemma 28. Then we can replace any $P_m + P_6 wr S_t$ by $P_{m+6} + P_6 wr S_{t-1}$ and similarly any $P_m + P_r$ by a symmetric group for $3 \leq r \leq 5$. □

**Proposition 31.** *Theorem A$'$ holds for $\ell = 2$.*

*Proof.* By Lemmas 29 and 30, if $n \geq 36$ we see that $G \cong S_{n+1}$ or $S_{n+2}$, and the bound holds. If no primitive constituent is a symmetric group, the bound fails for $n = 33$ since $(6531840)^5 \times 5! \times 216 = 3.08 \times 10^{38}$ and $34! = 2.95 \times 10^{38}$, but the claimed bounds do hold for $n = 34$ or $35$. □

## 9. The case $\ell = 5$

By Lemma 20, there are no primitive components of subdegrees 5, 7 or 9. The actual calculations are different from before since we have $[P_{2,5} : P_{2,5} Z(P_{2,5})] = 24$ and $[P_{3,5} : Z(P_{3,5})] = 2520$ so that in small degree primitive constituents of subdegree 3 dominate, rather than those of subdegree 2.

**Lemma 32.** *There is no primitive constituent of subdegree $m$ for $4 < m \leq 12$.*

*Proof.* The successive replacements of $P_{12}$ by $P_6 wr S_2$, $P_{11}$ by $P_{10} + P_1$, $P_{10}$ by $P_6 + P_4$, $P_8$ by $P_4 wr S_2$ and $P_6$ by $P_3 wr S_2$ would eliminate them. □

**Lemma 33.** (i) *Any primitive constituent that is a symmetric group has subdegree at least* 58.

(ii) *There are at most* 23 *primitive constituents of subdegree* 3.

(iii) *There is at most one primitive constituent of subdegree* 2. *If there is one, then there is none of subdegree* 1.

(iv) *There are at most five primitive constituents of subdegree* 4; *if there are five, there is none of subdegree* 1.

*Proof.* The inequalities $60! > (2520)^{19}.19! > 59!$ demonstrate the transition at which a constituent that is a symmetric group can no longer be replaced by $P_3 wr S_t$ for some $t$ to obtain a contradiction, while the reverse transition occurs when $70! < (2520)^{23}.23! < 71!$ This establishes (i) and (ii).

The replacement of $P_2 wr S_t$ by suitable $P_3 wr S_s + P_1 wr S_r$ if $2 \le t \le 12$, $P_2 wr S_t$ by $P_{2t}$ if $t > 12$, and $P_2 + P_1$ by $P_3$ establishes (iii), while replacing $P_4 wr S_t$ for $t > 5$ and $P_4 wr S_5 + P_1$ by suitable $P_3 wr S_{t'}$ yields (iv). $\square$

**Lemma 34.** *If any primitive constituent $P_m$ is a symmetric group, then there is no other nontrivial constituent.*

*Proof.* The inequality $60^3 > 5.25920$ enables the replacement of $P_m + P_4 wr S_t$ by $P_{m+4} + P_4 wr S_{t-1}$ to eliminate constituents of subdegree 4.

A similar, but more delicate, argument eliminates primitive constituents of subdegree 3; in the extreme case, we must use the inequality

$$61 \cdot 62 \cdots 68 > 2520^3 \cdot 23 \cdot 22 \cdot 21$$

to remove three constituents of degree 3 simultaneously if we had started with twenty three, replacing $P_{58} + P_3 wr S_{23}$ by $P_{67} + P_3 wr S_{20}$. We replace $P_m + P_2$ by $P_{m+2}$ to eliminate $P_2$. $\square$

**Proposition 35.** *Theorem $A'$ holds for $\ell = 5$.*

*Proof.* Since $2520^{23} \cdot (23)! > 70!$, we may suppose that $n \geq 70$. If any primitive constituent is a symmetric group, then Lemma 34 gives the desired conclusion. If not, we can replace any number of primitive constituents of subdegree 4 by additional constituents of subdegree 3; then Lemma 33 forces $n \leq 72$ and we can eliminate these possibilities. $\square$

## 10. Small degrees

In view of Theorem A, the determination of the best possible bounds for degrees less than 71 (or for $n < n(\ell)$ when $\ell \leq 5$ by Theorem A$'$) is the strictly finite problem of comparing all saturated pairs in which the primitive constituents appear in the conclusion of Theorem B. Consequently, since this search could be carried out using a computer, we will just indicate the necessary comparisons that eliminate such candidates in the same way as we carried out reductions in determining the bound in Theorem A, and we shall state the bounds in terms of groups that achieve them.

We refine some of our notation for the remainder of this paper; for a primitive pair $(P_m, k^m)$, we write $P_m^{(t)}$ to denote either the group $P_m wr S_t$ or the pair $(P_m wr S_t, k^{mt})$, according to context, and use $+$ both in its previous role for pairs, and to denote the direct product of groups. Let $f(n, \ell)$ denote the optimal bound that should replace that in Theorem A for $n < n(\ell)$ (or $n < 71$).

**Case I - $\ell \geq 7$.**

Since the primitive groups that arise in Theorem B when $\ell \geq 7$ are precisely those that occur in characteristic 0 with the exception of getting the symmetric group $S_{n+2}$ when $n \equiv -2 \pmod{\ell}$, we need only modify the characteristic 0 bounds obtained in [**4**] to allow for this. However, by Lemma 23, no symmetric group arises as a primitive constituent unless its subdegree is at least 53. Thus we can confine our attention to the range $53 \leq n \leq 70$ when either $G$ is a symmetric group, or else $G \cong (2 \cdot A_5) wr S_t$, allowing the further possibility that $G = P_n \cong S_{n+2}$ if $n \equiv -2 \pmod{\ell}$. This results in the following bounds.

**Theorem C.** *Let $k$ is an algebraically closed field of characteristic $\ell \geq 7$. Suppose that $n \leq 70$. Then $f(n,\ell)$ is given by one of the following, where $G$ is a finite subgroup of $GL(n,k)$ that achieves the bound $f(n,\ell)$.*

(i) $\quad n \leq 6$ *and $G$ is primitive.*

(ii) $\quad 7 \leq n \leq 19$, *$G$ is imprimitive, and the following saturated pairs achieve the bound $f(n,\ell)$:*

| $n$ | | $n$ | |
|---|---|---|---|
| 7 | $P_4 + P_3$ | 14 | $P_2^{(7)}$ |
| 8 | $P_4^{(2)}$ | 15 | $P_4^{(3)} + P_3$ |
| 9 | $P_6 + P_3$ | 16 | $P_4^{(4)}$ |
| 10 | $P_6 + P_4$ | 17 | $P_4^{(4)} + P_1$ |
| 11 | $P_4^{(2)} + P_3$ | 18 | $P_2^{(9)}$ |
| 12 | $P_4^{(3)}$ | 19 | $P_4^{(4)} + P_3$ |
| 13 | $P_4^{(3)} + P_1$ | | |

(iii) *If $20 \leq n \leq 70$ and $n = 2r$ or $2r+1$, the bounds $f(n,\ell)$ are as follows:*

| $f(n,\ell)$ | $n$ | *restriction on $\ell$* |
|---|---|---|
| $(n+2)!$ | $53 \leq n \leq 60$, $n$ odd | $\ell \mid n+2$ |
| $(n+2)!$ | $61 \leq n \leq 70$ | $\ell \mid n+2$ |
| $(n+1)!$ | $n \geq 63$, $n$ odd | $\ell \nmid n+2$ |
| $60^r \cdot r!$ | *otherwise* | |

*These bounds are achieved by $P_n$, $P_{n-1} + P_1$, $P_2^{(r)}$ or $P_2^{(r)} + P_1$.*

**Case II - $\ell = 2$.**

By Theorem $A'$, we can confine our attention to the range $n \leq 33$. By Lemma 20 and the reductions of Section 8, we may suppose that that either there is just one nontrivial primitive constituent that is a symmetric group, or else the bound is given by a group whose shape is $P_{6,2}^{(r)}$, $P_{6,2}^{(r)} + P_{d,2}$ for $d = 1, 3, 4$ or $5$, or $P_{6,2}^{(r)} + (P_1, k) wr S_2$, and then $r \leq 5$. It is now a simple calculation to find the bounds, with $[P_{3,2} : Z(P_{3,2})] = 216$, $[P_{4,2} : Z(P_{4,2})] = 2520$, $[P_{5,2} : Z(P_{5,2})] = 3000$ and $[P_{6,2} : Z(P_{6,2})] = 6531840$. We note also that we have the "exceptional" primitive constituent $S_{n+2}$ in degree $n$ only for $n$ even.

**Theorem D.** *Suppose that $n = 6r + d \leq 33$ where $0 \leq d \leq 5$. Then $f(n, 2)$ is given by saturated pairs of the form $P_{6,2}^{(r)}$, $P_{6,2}^{(r)} + P_{d,2}$ for $d = 1, 3, 4$ or $5$, or $P_{6,2}^{(r)} + (P_1, k) wr S_2$, unless $n = 26, 28, 29$ or $32$; $f(n, 2) = (n+1)!$ or $(n+2)!$ in the exceptional cases, as $n$ is odd or even, respectively.*

**Case III - $\ell = 3$.**

By Lemma 26, any primitive constituent not a symmetric group has subdegree at most 4, while at most one is a symmetric group and then of subdegree at least 55 by Lemma 22 and the argument of Proposition 28. We need consider only the cases when $n \leq 68$. Now $[P_{4,3} : Z(P_{4,3})] = 40320$, $[P_{3,3} : Z(P_{3,3})] = 168$ and $[P_{2,3} : Z(P_{2,3})] = 60$; a consequence is that the interaction between primitive constituents of degrees 2 and 4 become far more delicate than in the characteristic 0 case.

The replacements $P_3^{(r)}$ by $P_2^{(s)}$ if $r \geq 2$, $P_2^{(t)} + P_3$ by $P_2^{(t+1)} + P_1$ if $t \geq 2$ and $P_2 + P_3$ by $P_4 + P_1$ show that there is at most one primitive constituent of subdegree 3, and none if there is any of subdegree 2 while the inequalities $60^{2t}.(2t)! > 40320^t.t!$ if $t \geq 8$ and $60^{14}.(4)! < 40320^7.7!$ show that there are at most seven of subdegree 4. A series of comparisons now shows that if any has subdegree 4, then there is at most one of subdegree 2. Further comparision with the characteristic 0 case for $n \geq 55$ when $n \equiv 1 \pmod 3$ yields the following conclusion.

**Theorem E.** *Suppose that $n \leq 68$. Then the bound $f(n, 3)$ is given by the saturated pairs in the following table:*

| $n$ | exclusions | saturated pairs |
|---|---|---|
| $n \leq 4$ | | $P_{n,3}$ |
| $5 \leq n = 4r + d \leq 29 \ (d \leq 3)$ | $n = 22, 26, 27$ | $P_{4,3}^{(r)}, \ P_{4,3}^{(r)} + P_{d,3}$ |
| $55 \leq n = 3m + 1 \leq 67$ | | $P_{n,3} \ (\cong S_{n+2})$ |
| $57$ | | $P_{55,3} + P_{2,3} \ (\cong S_{56} \times 2 \cdot A_5)$ |
| $63, \ 65$ | | $P_{n,3} \ (\cong S_{n+1})$ |
| $n = 2r \ or \ 2r+1, \ otherwise$ | | $P_{2,3}^{(r)} \ or \ P_{2,3}^{(r)} + P_1$ |

**Case IV - $\ell = 5$.**

By Lemma 32, every primitive constituent that is not a symmetric group has subdegree at most 4; hence, as the principles of the calculations are unaltered and relatively few comparisons need be made, we omit the details. Here,

$$[P_{2,5} : Z(P_{2,5})] = 24, \ [P_{3,5} : Z(P_{3,5})] = 2520 \text{ and } [P_{4,5} : Z(P_{4,5})] = 25920 .$$

**Theorem F.** *Suppose that $n \leq 69$. Then the bound $f(n,5)$ is given by the saturated pairs in the following table*:

| $n$ | exclusions | saturated pairs |
|---|---|---|
| $n \leq 4$ | | $P_{n,5}$ |
| $n = 5$ | | $P_{3,5} + P_{2,5}$ |
| $n = 8$ | | $P_{4,5}^{(2)}$ |
| $n = 11$ | | $P_{4,5}^{(2)} + P_{3,5}$ |
| $6 \leq n = 3r \leq 69$ | | $P_{3,5}^{(r)}$ |
| $7 \leq n = 3r + 1 \leq 31$ | | $P_{3,5}^{(r-1)} + P_{4,5}$ |
| $34 \leq n = 3r + 1 \leq 64$ | $n = 58$ | $P_{3,5}^{(r)} + P_1$ |
| $14 \leq n = 3r + 2 \leq 62$ | | $P_{3,5}^{(r)} + P_{2,5}$ |
| $n = 65, \ 67$ | | $P_{n,5} \ (\cong S_{n+1})$ |
| $n = 58, \ 68$ | | $P_{n,5} \ (\cong S_{n+2})$ |

## 11. Weisfeiler's bounds

Had we solely determined numerical orders for the bounds given in Theorems C – F by use of a computer, the interplay between primitive constituents of the various subdegrees would have been be lost, and in particular we would not see the reasons for the extent to which the bounds are by no means "smooth", even for a fixed field.

The bound of $n^4(n+2)!$ that Weisfeiler gave for $n < 64$ in his manuscript is very much a smooth overarching bound and reflects the concept behind his earlier, weaker, bounds in [**15**] where he separated out the contributions from components that were of alternating type from those that were of Lie type. It actually gets closest when $n = 18$ and $\ell = 5$ where

$$18^4 \cdot 20! = 2.55 \times 10^{23} \text{ and } (2520)^6.6! = 1.84 \times 10^{23}$$

when a "best" bound would be given by $n^\alpha \times (n+2)!$ with $\alpha = 3.8873$. However, this value of $\alpha$ is far from optimal in characteristic $\ell \neq 5$ or for different values of $n$; for example, if $\ell \neq 5$ and $n \geq 60$ (or $n \geq 34$ if $\ell = 2$) we could take $\alpha = 0$.

University College
Oxford OX1 4BH
England

e-mail:  mjc@herald.ox.ac.uk